\documentclass[journal]{IEEEtran}
\usepackage[utf8]{inputenc}
\usepackage[english]{babel}

\usepackage[colorlinks=true, allcolors=blue]{hyperref}

\usepackage{amsmath, amssymb}

\usepackage{algorithm, algpseudocode}

\newtheorem{theorem}{Theorem}
\newtheorem{lemma}{Lemma}

\newtheorem{remark}{Remark}
\newtheorem{definition}{Definition}
\newtheorem{proposition}{Proposition}

\usepackage{graphicx}
\graphicspath{ {figure/} }
\usepackage{subfigure}
\usepackage{tabularx}
\usepackage{threeparttable}
\usepackage{booktabs} 
\usepackage{multirow}
\usepackage[noadjust]{cite}
\usepackage{url}

\IEEEoverridecommandlockouts

\title{\LARGE \bf
A Gap Penalty Reformulation for Mathematical Programming with Complementarity Constraints: Convergence Analysis}
\author{Kangyu Lin and Toshiyuki Ohtsuka
\thanks{This work was supported by the JSPS KAKENHI (No. JP22H01510 and JP23K22780). 
Kangyu Lin was supported by the CSC scholarship (No. 201906150138).
The authors are with the Graduate School of Informatics, Kyoto University, Japan.
Corresponding author: Kangyu Lin
(Email: {\tt\small k-lin@sys.i.kyoto-u.ac.jp}.)}
}

\begin{document}

\maketitle

\begin{abstract}
Our recent study \textcolor{blue}{(Lin and Ohtsuka, 2024)} proposed a new penalty method for solving mathematical programming with complementarity constraints (MPCC).
This method first reformulates MPCC as a parameterized nonlinear programming called gap penalty reformulation and then solves a sequence of gap penalty reformulations with an increasing penalty parameter.
This study examines the convergence behavior of the new penalty method.
We prove that it converges to a strongly stationary point of MPCC, provided that: 
(i) The MPCC linear independence constraint qualification holds. 
(ii) The upper-level strict complementarity condition holds. 
(iii) The gap penalty reformulation satisfies the second-order necessary conditions in terms of the second-order directional derivative.
Because strong stationarity is used to identify the MPCC local minimum, our analysis indicates that the new penalty method can find an MPCC solution.


\end{abstract}

\section{Introduction}

\subsection{Background}

This study considers a special yet common class of nonlinear programming (NLP), known as \textit{mathematical programming with complementarity constraints} (MPCC).
\textcolor{blue}{Many switching and nonsmooth decisions can be represented through complementarity constraints.
Therefore, MPCC can model the discretized optimal control problems of various nonsmooth dynamical systems arising in practical applications, such as process systems with discrete events, mechanical systems with contacts, and game-theoretical dynamical systems \cite{biegler2010nonlinear}}.

However, complementarity constraints pose significant challenges in solving the MPCC, as almost all constraint qualifications (CQs) are violated at any MPCC feasible point.
This causes two difficulties.
First, MPCC solutions cannot be characterized by the optimality conditions for standard NLP problems;
Second, MPCC cannot be solved using standard NLP solution methods.
To characterize the MPCC solution, many MPCC-tailored concepts have been proposed \cite{scholtes2001convergence}, some of which are reviewed in Section \ref{subsection: MPCC-tailored concepts}.
\textcolor{blue}{To solve the MPCC, many MPCC-tailored solution methods, such as \textit{relaxation methods} and \textit{penalty methods}, have been proposed \cite{biegler2010nonlinear}}.
These methods do not solve the MPCC directly. 
Instead, they reformulate complementarity constraints as parameterized costs or constraints, solve a sequence of well-defined parameterized NLP, and then present a rigorous convergence analysis showing that the sequence of solutions to the parameterized NLP can finally converge to an MPCC solution.
These MPCC-tailored solution methods are practical because they are easy to implement using state-of-the-art NLP software.

Our recent study \cite{lin2024gap} proposed a new penalty method called \textit{gap penalty method}, which reformulates the MPCC as a parameterized NLP called \textit{gap penalty reformulation} and solves a sequence of gap penalty reformulations with an increasing penalty parameter.
\textcolor{blue}{Compared to other penalty methods, gap penalty reformulation exhibits certain convexity structures that can be exploited by a dedicated Hessian regularization method.}
\textcolor{blue}{The gap penalty method is practically effective, especially in optimal control of linear complementarity systems} \cite{lin2024gap}, but its theoretical convergence has not been analyzed.

\subsection{Contribution}

This study builds on our previous work \cite{lin2024gap} by providing a rigorous convergence analysis for the gap penalty method.
This study’s main contributions are:
First, because the gap penalty reformulation involves a \textit{Lipschitz continuously differentiable} function, we characterize its solutions using \textit{second-order directional derivatives} rather than the Hessian matrix as in standard NLP theory.
These results are critical for the subsequent convergence analysis.
Second, we prove that the gap penalty method converges to a \textit{strongly stationary point} of MPCC under certain standard MPCC-tailored assumptions. 
This \textit{theoretically} confirms that the gap penalty method can find an MPCC solution.

\subsection{Outline and notation}

The remainder of this study is organized as follows:
Section \ref{section: mathematical programming with complementarity constraints} reviews MPCC-tailored concepts and gap penalty reformulation;
Section \ref{section: lipschitz continuously differentiable Optimization} presents the optimality conditions of gap penalty reformulation;
Section \ref{section: convergence analysis} presents a detailed convergence analysis;
and Section \ref{section: conclusion} concludes this study.

\textcolor{blue}{Given a vector $x \in \mathbb{R}^n$, we denote its $i$-th component by $x_i$}, \textcolor{blue}{its $\ell_1$ and $\ell_2$ norm by $\| x\|_1$ and $ \| x\|$, respectively}.
We denote the complementarity condition between $x, y \in \mathbb{R}^n$ as $0 \leq x \perp y \geq 0$, where $x \perp y$ means $x^T y = 0$.
\textcolor{blue}{We take $\max(x, y)$ with $x, y \in \mathbb{R}^n$ in a component-wise manner}.
\textcolor{blue}{Given a function $f: \mathbb{R}^n \rightarrow \mathbb{R}^m$, we denote its Jacobian by $\nabla_x f \in \mathbb{R}^{m \times n}$}
and say that $f$ is $k$-th Lipschitz continuously differentiable ($LC^k$) if its $k$-th derivative is Lipschitz continuous.

\section{Mathematical programming with complementarity constraints}\label{section: mathematical programming with complementarity constraints}
\subsection{MPCC-tailored concepts}\label{subsection: MPCC-tailored concepts}
Consider the MPCC in the form of:
\begin{subequations}\label{equation: MPCC}
    \begin{align}
        \min_{x, \lambda, \eta} \ & J(x, \lambda, \eta), \label{equation: MPCC cost function}\\
        \text{s.t.} \ & h(x, \lambda, \eta) = 0, \label{equation: MPCC equality constraints}\\
                      & 0 \leq \lambda \perp \eta \geq 0, \label{equation: MPCC complementarity constraints} 
    \end{align}
\end{subequations}
where $x \in \mathbb{R}^{n_x}$, $\lambda, \eta \in \mathbb{R}^{n_{\lambda}}$ are decision variables,
$J: \mathbb{R}^{n_x} \times \mathbb{R}^{n_{\lambda}} \times \mathbb{R}^{n_{\lambda}} \rightarrow \mathbb{R}$ is the cost function, 
and $h: \mathbb{R}^{n_x} \times \mathbb{R}^{n_{\lambda}} \times \mathbb{R}^{n_{\lambda}} \rightarrow \mathbb{R}^{n_h}$ is the equality constraint. 
$J$ and $h$ are assumed to be $LC^2$.
We define $z = [x^T, \lambda^T, \eta^T]^T \in \mathbb{R}^{n_x + 2 n_{\lambda}}$ to collect decision variables.
A point $z$ satisfying (\ref{equation: MPCC equality constraints}) and (\ref{equation: MPCC complementarity constraints}) is referred to as the \textit{feasible point} of MPCC (\ref{equation: MPCC}).
For a feasible point $z$ of MPCC (\ref{equation: MPCC}), the activation status of (\ref{equation: MPCC complementarity constraints}) at $z$ are classified by the following index sets:
\begin{subequations}\label{equation: MPCC index set}
    \begin{align}
       & \mathcal{I}_{\lambda}(z) = \{ {i \in \{1, \cdots, n_{\lambda} \} | \lambda_i = 0, \eta_i > 0 \}}, \\
       & \mathcal{I}_{\eta}(z) = \{ {i \in \{1, \cdots, n_{\lambda} \} | \lambda_i > 0, \eta_i = 0 \}}, \\
       & \mathcal{I}_{\lambda \eta}(z) = \{ {i \in \{1, \cdots, n_{\lambda} \} | \lambda_i = 0, \eta_i = 0 \}}.
    \end{align}
\end{subequations}
Subsequently, the MPCC-tailored stationarity, constraint qualification, and strict complementarity conditions can be defined based on these index sets \cite{scholtes2001convergence}:
\begin{definition}[MPCC-tailored stationarity]
    For a feasible point $z$ of MPCC (\ref{equation: MPCC}), we say that 
    \begin{itemize}
        \item it is \textit{weakly stationary} if there exists Lagrangian multipliers $u \in\mathbb{R}^{n_h}$, $v,w \in \mathbb{R}^{n_{\lambda}}$ for $h(z) = 0$, $\lambda \geq 0$ and $\eta \geq 0$, respectively, such that:
            \begin{subequations}\label{equation: weak stationary definition}
                \begin{align}
                    & \nabla_z J(z) + u^T \nabla_z h(z) - v^T \nabla_z \lambda - w^T \nabla_z \eta = 0, \\
                    & v_i \in \mathbb{R}, w_i = 0,                           \ i \in \mathcal{I}_{\lambda}(z), \\
                    & v_i = 0,            w_i \in \mathbb{R},  \ i \in \mathcal{I}_{\eta}(z), \\
                    & v_i \in \mathbb{R}, w_i \in \mathbb{R},  \ i \in \mathcal{I}_{\lambda \eta} (z),
                \end{align}
            \end{subequations}
        \item it is \textit{Clarke stationary} if it is weakly stationary and $ v_i w_i \geq 0, i \in \mathcal{I}_{\lambda \eta} (z)$,
        \item it is \textit{strongly stationary} if it is weakly stationary and $v_i\geq 0, w_i \geq 0, i \in \mathcal{I}_{\lambda \eta} (z)$.
    \end{itemize}
\end{definition}

\begin{definition}[MPCC-tailored constraint qualification]
    For a feasible point $z$ of MPCC (\ref{equation: MPCC}), we say that the \textit{MPCC linear independence constraint qualification (MPCC-LICQ)} holds at $z$ if the vectors $\{\nabla_z h (z)  \} \cup \{ \nabla_z \lambda_i | i \in  \mathcal{I}_{\lambda} (z) \cup \mathcal{I}_{\lambda \eta}  (z) \} \cup \{ \nabla_z \eta_i | \mathcal{I}_{\eta} (z) \cup \mathcal{I}_{\lambda \eta} (z)\}$ are linearly independent.
\end{definition}

\begin{definition}[MPCC-tailored strict complementarity condition]
    A weakly stationary point $z$ is said to satisfy the \textit{upper-level strict complementarity (ULSC) condition} if there exist multipliers satisfying (\ref{equation: weak stationary definition}) and $ v_i w_i \neq 0,  \forall i \in \mathcal{I}_{\lambda \eta} (z)$.
\end{definition}

With these MPCC-tailored concepts, the solution of MPCC (\ref{equation: MPCC}) can be characterized by the following \textit{first-order necessary optimality condition for MPCC}. 
\begin{proposition}[Theorem 11.1 in \cite{biegler2010nonlinear}]
    Let $\bar{z}$ be a local minimum of MPCC (\ref{equation: MPCC}) and suppose that the MPCC-LICQ holds at $\bar{z}$, then $\bar{z}$ is a strongly stationary point of MPCC (\ref{equation: MPCC}).
\end{proposition}

\subsection{Gap penalty reformulation for MPCC}
The treatment of complementarity constraints plays a critical role in MPCC-tailored solution methods.
Our recent study \cite{lin2024gap} uses the D-gap function \cite{peng1997unconstrained} to construct a novel penalty term for the complementarity constraints (\ref{equation: MPCC complementarity constraints}).
\begin{definition}[D-gap function]
    Let $\lambda, \eta \in \mathbb{R}^{n_{\lambda}}$ be two variables, $a, b$ be two given constants satisfying $b > a > 0$, and $\varphi^{ab}: \mathbb{R}^{n_{\lambda}} \times \mathbb{R}^{n_{\lambda}} \rightarrow \mathbb{R}$ be a function given by
    \begin{equation}\label{equation: definition of the D gap function}
        \varphi^{ab}(\lambda, \eta) = \varphi^{a}(\lambda, \eta) - \varphi^{b}(\lambda, \eta),
    \end{equation}
    where $\varphi^{a}, \varphi^{b}: \mathbb{R}^{n_{\lambda}} \times \mathbb{R}^{n_{\lambda}} \rightarrow \mathbb{R}$ are the functions defined by
    \begin{equation}\label{equation: elements of the D gap function}
        \varphi^{c}(\lambda, \eta) = \frac{1}{2 c} (\| \eta \|^2_2 - \| \max(0, \eta -c \lambda) \|^2_2),
    \end{equation}
    with parameter $c > 0$.
    We call $\varphi^{ab}$ the \textit{D-gap function}, where D stands for Difference.
\end{definition}

Some properties of $\varphi^{ab}$ are summarized below;
see Theorem 10.3.3 and Proposition 10.3.13 in \cite{facchinei2003finite} for a detailed proof.
\begin{proposition}\label{proposition: properties of D gap function}
    The following statements are valid for $\varphi^{ab}$.
    \begin{itemize}
        \item (\textit{Differentiability}) $\varphi^{ab}(\lambda, \eta)$ is $LC^1$ with the gradients
        \begin{subequations}\label{equation: explicit expression of D gap gradient}
            \begin{align}
                \begin{split}
                    \nabla_{\lambda}\varphi^{ab} = & [\max(0, \eta - a \lambda)  - \max(0, \eta - b \lambda)]^T
                \end{split}     
                \\
                \begin{split}
                    \nabla_{\eta}\varphi^{ab} = & [(\frac{1}{a} - \frac{1}{b})\eta - \frac{1}{a} \max (0, \eta - a \lambda)  \\ 
                                                            & + \frac{1}{b} \max(0, \eta - b\lambda)]^T. 
                \end{split}                              
            \end{align}
        \end{subequations}    
        \item (\textit{Non-negativity}) $\varphi^{ab}(\lambda, \eta) \geq 0, \forall \lambda, \eta \in \mathbb{R}^{n_{\lambda}}$.
        \item (\textit{Equivalence}) 
        $\varphi^{ab}(\lambda, \eta) = 0 $ iff $0 \leq \lambda \perp \eta \geq 0$.
    \end{itemize}
\end{proposition}

Note that $\varphi^{ab}$ exhibits \textit{partial separability} (Definition 7.1, \cite{gros2020numerical}), and can be written as the sum of $n_{\lambda}$ scalar subfunctions:
\begin{equation}\label{equation: reformulate D gap function as the sum of scalar functions}
    \varphi^{ab}(\lambda, \eta) = \sum^{n_{\lambda}}_{i = 1} \delta^{ab}(\lambda_{i}, \eta_{i}),
\end{equation}
where subfunction $\delta^{ab}: \mathbb{R} \times \mathbb{R} \rightarrow \mathbb{R}$ is given by:
\begin{equation}\label{equation: scalar D gap function explicit expression}
    \begin{split}
         \delta^{ab}(\lambda_{i}, \eta_{i}) = & \frac{b-a}{2ab}\eta_{i}^2 - \frac{1}{2a}\{\max(0, \eta_{i} - a\lambda_{i})\}^2 \\
                                            & + \frac{1}{2b} \{\max(0, \eta_{i} - b\lambda_{i})\}^2,
    \end{split}    
\end{equation}
which is a special case of $\varphi^{ab}$ with variables being scalars.

By replacing the complementarity constraints (\ref{equation: MPCC complementarity constraints}) with the D-gap function, we obtain a parameterized NLP problem \cite{lin2024gap}
\begin{subequations}\label{equation: MPCC gap penalty reformulation}
    \begin{align}
        \mathcal{P}_{gap}(\mu): \quad \quad \min_{x, \lambda, \eta} \ & J(x, \lambda, \eta) + \mu\varphi^{ab}(\lambda, \eta), \label{equation: MPCC gap penalty reformulation cost function}\\
        \text{s.t.} \ & h(x, \lambda, \eta) = 0. \label{equation: MPCC gap penalty reformulation equality constraints} 
    \end{align}
\end{subequations}
with penalty parameter $\mu > 0$.
$\mathcal{P}_{gap}(\mu)$ is referred to as the \textit{gap penalty reformulation} for MPCC (\ref{equation: MPCC}).
We hope that the solution to MPCC (\ref{equation: MPCC}) can be obtained by solving a sequence of $\mathcal{P}_{gap}(\mu)$ with $\mu \rightarrow + \infty$.
However, since $\mathcal{P}_{gap}(\mu)$ involves an $LC^1$ function $\varphi^{ab}$, its solution cannot be identified using the Hessian matrix, as in standard NLP theory.
Therefore, we must discuss how to determine the solution to $\mathcal{P}_{gap}(\mu)$.
These discussions constitute the first contribution of this paper and are essential for the subsequent convergence analysis.

\section{Lipschitz continuously differentiable Optimization}\label{section: lipschitz continuously differentiable Optimization}

In this section, we use the directional derivatives to identify the solution to an NLP problem with $LC^1$ functions.

\subsection{Directional derivatives}

First, we provide the definitions and properties of the first- and second-order directional derivatives for a function \cite{ben1982necessary}.

\begin{definition}[Directional derivatives]
    Let $f:  \mathbb{R}^{n} \rightarrow  \mathbb{R}$ be a real-valued function.
    \begin{itemize}
        \item The \textit{first-order directional derivative} of $f$ at $x$ in the direction $d$ is defined as (if this limit exists)
        \begin{equation}
            D(f(x); d) = \lim_{t \downarrow 0} \frac{f(x + td) - f(x)}{t}.
        \end{equation}       
        \item Supposing that $D(f(x); d)$ exists, the \textit{second-order directional derivative} of $f$ at $x$ in the directions $d$ and $p$ is defined as (if this limit exists)
        \begin{equation}
            \begin{split}
                 & D^2(f(x); d, p) \\
                = &  \lim_{t \downarrow 0} \frac{f(x + td + t^2p) - f(x) - t D(f(x); d)}{t^2}.
            \end{split}           
        \end{equation}
    \end{itemize}
\end{definition}

\begin{proposition}\label{proposition: properties of the directional derivatives}
    Let $f:  \mathbb{R}^{n} \rightarrow  \mathbb{R}$ be a real-valued function.
    \begin{itemize}
        \item If $f$ is $C^1$, then $D(f(x); d)$ exists, and we have
        \begin{equation}
            D(f(x); d) = \nabla_x f(x) d.
        \end{equation}
        \item If $f$ is $C^2$, then $D^2(f(x); d, p)$ exists, and we have 
        \begin{equation}
            D^2(f(x); d, p) = \nabla_x f(x) p + \frac{1}{2}d^T \nabla_{xx} f(x) d.
        \end{equation}
        \item If $D(f(x); d)$ and $D^2(f(x); d, p)$ exist, then $f$ can be expanded in terms of directional derivatives with $t \geq 0$
        \begin{equation}\label{equation: function expansion in terms of directional derivatives}
            \begin{split}
                & f(x + td + t^2p) \\ 
                = & f(x) + tD(f(x); d) + t^2 D^2(f(x); d, p) + o(t^2).
            \end{split}          
        \end{equation}
        \item Let $f(x) = (\max(g(x), 0))^2$, where $g: \mathbb{R}^{n} \rightarrow  \mathbb{R}$ is affine.
        Then, $f$ is $LC^1$, and its second-order directional derivative at $x$ in the direction $d$ can be explicitly written as 
        \begin{equation}
            \begin{split}
                & D^2(f(x); d, d) \\
                = & 
                \begin{cases}
                    2 g(x) \nabla_x g(x)d + (\nabla_x g(x)d)^2, & \text{if } g(x) > 0 \\
                    (\max(0, \nabla_x g(x) d))^2,  & \text{if } g(x) = 0 \\
                    0,                     & \text{if } g(x) < 0 
                \end{cases}                
            \end{split}
        \end{equation}      
    \end{itemize}
\end{proposition}
\begin{IEEEproof}
    The first three statements are from Section 2 in \cite{ben1982necessary}, and the fourth is from Proposition 3.3 in \cite{ben1982necessary}.
\end{IEEEproof}

\subsection{Optimality condition}
Next, we discuss the first- and second-order necessary optimality conditions for an NLP problem with $LC^1$ functions.
Considering the NLP problem in the form of
\begin{subequations}\label{equation: NLP problem with LC1 function}
    \begin{align}
        \min_{x}    \ & J(x) + g(x), \label{equation: NLP problem with LC1 function cost function}\\
        \text{s.t.} \ & h(x) = 0, \label{equation: NLP problem with LC1 function equality constraints} 
    \end{align}
\end{subequations}
where
functions $J: \mathbb{R}^{n} \rightarrow \mathbb{R}$ and $h: \mathbb{R}^{n} \rightarrow \mathbb{R}^m$ are $LC^2$, 
and $g: \mathbb{R}^{n} \rightarrow \mathbb{R}$ is $LC^1$. 
A point $x$ satisfying (\ref{equation: NLP problem with LC1 function equality constraints}) is referred to as a feasible point of NLP (\ref{equation: NLP problem with LC1 function}).
Let $\gamma \in \mathbb{R}^m$ be the Lagrangian multiplier for (\ref{equation: NLP problem with LC1 function equality constraints}). 
The Lagrangian of (\ref{equation: NLP problem with LC1 function}) is defined as
\begin{equation}
    \mathcal{L}(x, \gamma) = J(x) + g(x) + \gamma^T h(x).
\end{equation}
The critical cone at a feasible point $x$ of (\ref{equation: NLP problem with LC1 function}) is defined as 
\begin{equation}
    \mathcal{C}(x) = \{ d \in \mathbb{R}^{n} | \nabla_x h(x) d = 0\}.
\end{equation}
We state the necessary optimality conditions of (\ref{equation: NLP problem with LC1 function}) below.
\begin{theorem}\label{theorem: necessary optimality conditions of the NLP problem with LC1 function}
    Let $x^*$ be a local minimum of NLP (\ref{equation: NLP problem with LC1 function}). 
    \begin{itemize}
        \item Assume that $\nabla_x h(x^*)$ has a full row rank. 
        Then the first-order necessary optimality condition of (\ref{equation: NLP problem with LC1 function}) holds at $x^*$, that is, there exists a unique multiplier $\gamma^*$ such that 
        \begin{subequations}\label{equation: NLP problem with LC1 function first order necessary optimality condition}
            \begin{align}
                & \underbrace{\nabla_x J(x^*) + \nabla_x g(x^*) + (\gamma^*)^T \nabla_x h(x^*)}_{\nabla_x\mathcal{L}(x^*, \gamma^*)} = 0, \\
                & h(x^*) = 0.
            \end{align}
        \end{subequations} 
        \item Furthermore, assume that $D^2(g(x^*); d, d)$ exists for any $d \in \mathcal{C}(x^*)$. 
        Then the second-order necessary optimality condition of (\ref{equation: NLP problem with LC1 function}) in terms of the second-order directional derivatives holds at $x^*$, that is, 
        \begin{equation}\label{equation: NLP problem with LC1 function second order necessary optimality condition}
            D^2(J(x^*) + g(x^*); d, d) \geq 0, \textcolor{blue}{\forall  d \in \mathcal{C}(x^*)}
        \end{equation}
    \end{itemize}
\end{theorem}
\begin{IEEEproof}
    The conditions (\ref{equation: NLP problem with LC1 function first order necessary optimality condition}) are also known as the \textit{Karush-Kuhn-Tucker (KKT) conditions}, where the proof can be found in many monographs (e.g., Section 12.4 in \cite{nocedal2006numerical}). 
    The second statement is proved in Appendix \ref{subsection: proof of necessary optimality conditions of the NLP problem with LC1 function}. 
\end{IEEEproof}

\begin{remark}
    The necessary optimality condition for \textit{unconstrained} NLP problems with $LC^1$ functions was derived in \cite{ben1982necessary}, 
    but we did not find any related derivations for the \textit{constrained} case (\ref{equation: NLP problem with LC1 function}). 
    Thus, Theorem \ref{theorem: necessary optimality conditions of the NLP problem with LC1 function} extends the work in \cite{ben1982necessary}.  
\end{remark}

Finally, a pair $(x, \gamma)$ satisfying (\ref{equation: NLP problem with LC1 function first order necessary optimality condition}) is referred to as a \textit{KKT point} of NLP (\ref{equation: NLP problem with LC1 function}), and the primal part $x$ is referred to as a \textit{stationary point} of NLP (\ref{equation: NLP problem with LC1 function}).

\section{Convergence analysis}\label{section: convergence analysis}

Let $z^k$ be a stationary point of $\mathcal{P}_{gap}(\mu^k)$ with given $\mu^k > 0$.
In this section, we analyze the convergence behavior of the sequence $\{z^k\}_{k = 0}^{\infty}$.
Let $\bar{z}$ be a limit point of $\{z^k\}_{k = 0}^{\infty}$ with $\mu^k \rightarrow + \infty$. 
Through the following three steps, each involving progressively stronger assumptions, we prove that $\bar{z}$ is a strongly stationary point of MPCC (\ref{equation: MPCC}).

\textit{Step 1}: prove that $\bar{z}$ is a feasible point of MPCC;

\textit{Step 2}: prove that $\bar{z}$ is a Clarke stationary point of MPCC;

\textit{Step 3}: prove that $\bar{z}$ is a strongly stationary point of MPCC.

\subsection{Feasibility analysis}

We begin by discussing the feasibility of limit point $\bar{z}$.
\begin{lemma}\label{lemma: feasibility of limit point}
    For a given $\mu^k > 0$, let $z^k$ be a stationary point of $\mathcal{P}_{gap}(\mu^k)$.
    Let $\mu^k \rightarrow + \infty$ and $\bar{z}$ be a limit point of $\{z^k\}_{k = 0}^{\infty}$.
    Assume that the value of cost function (\ref{equation: MPCC gap penalty reformulation cost function}) associated with $\{z^k\}_{k = 0}^{\infty}$ is bounded, that is, there exists a real number $M$ such that $| J(z^k) + \mu^k\varphi^{ab}(z^k)| \leq M, \forall k$.
    Then, the limit point $\bar{z}$ is a feasible point of MPCC (\ref{equation: MPCC}).
\end{lemma}
\begin{IEEEproof}
    The proof is inspired by Theorem 4.2 in \cite{huang2006sequential}.
    Since $\bar{z}$ is a limit point of $\{z^k\}_{k = 0}^{\infty}$, there exist a subsequence $\mathcal{K}$ such that $\lim\limits_{k \rightarrow \infty} z^k = \bar{z}$ for $k \in \mathcal{K}$.
    It is clear that $\bar{z}$ satisfies $h(\bar{z}) = 0$ because $z^k$ is the solution of (\ref{equation: MPCC gap penalty reformulation}). 
    Therefore, we only need to show that $\bar{z}$ satisfies the complementarity constraints $0 \leq \bar{\lambda} \perp \bar{\eta} \geq 0$.
    Since $| J(z^k) + \mu^k\varphi^{ab}(z^k)| \leq M, \forall k$, there exists a real number $M_1$ such that
    \begin{equation}\label{equation: boundness of the gap penalty function}
        \varphi^{ab}(z^k) \leq \frac{M_1}{\mu^k}.
    \end{equation} 
    Since $\varphi^{ab}(z^k) \geq 0$, by taking the limit of (\ref{equation: boundness of the gap penalty function}) as $k \rightarrow \infty$ for $k \in \mathcal{K}$, we have $\varphi^{ab}(\bar{z}) = 0$.
    Thus, following from the third statement in Proposition \ref{proposition: properties of D gap function}, we have $0 \leq \bar{\lambda} \perp \bar{\eta} \geq 0$.
\end{IEEEproof}

\subsection{Clarke stationarity analysis}
Next, we show that the limit point $\bar{z}$ is a Clarke stationary point of MPCC (\ref{equation: MPCC}) if the MPCC-LICQ holds at $\bar{z}$.

\begin{theorem}\label{theorem: Clarke stationarity of limit point}
    For a given $\mu^k > 0$, let $z^k$ be a stationary point of $\mathcal{P}_{gap}(\mu^k)$.
    Let $\mu^k \rightarrow + \infty$ and $\bar{z}$ be a limit point of $\{z^k\}_{k = 0}^{\infty}$.
    Let the assumption of Lemma \ref{lemma: feasibility of limit point} hold, and assume that the MPCC-LICQ holds at $\bar{z}$.
    Then, $\bar{z}$ is a Clarke stationary point of the MPCC (\ref{equation: MPCC}), that is, there exist Lagrangian multipliers $\bar{u} \in\mathbb{R}^{n_h}$, $\bar{v},\bar{w} \in \mathbb{R}^{n_{\lambda}}$ such that:
        \begin{subequations}\label{equation: Clarke stationarity analysis equation}
            \begin{align}
                & \nabla_z J(\bar{z}) + \bar{u}^T \nabla_z h(\bar{z}) - \bar{v}^T \nabla_z \bar{\lambda} - \bar{w}^T \nabla_z \bar{\eta} = 0, \label{equation: Clarke stationarity analysis equation stationary condition}\\
                & \bar{v}_i \in \mathbb{R}, \bar{w}_i = 0,    \ i \in \mathcal{I}_{\lambda}(\bar{z}), \label{equation: Clarke stationarity analysis equation active lambda}\\
                & \bar{v}_i = 0,   \bar{w}_i \in \mathbb{R},  \ i \in \mathcal{I}_{\eta}(\bar{z}), \label{equation: Clarke stationarity analysis equation active eta}\\
                & \bar{v}_i \bar{w}_i \geq 0,                 \ i \in \mathcal{I}_{\lambda \eta}(\bar{z}) \label{equation: Clarke stationarity analysis equation active lambda eta}.
            \end{align}
        \end{subequations}    
\end{theorem}

\begin{IEEEproof}
    The proof is partially inspired by Lemma 4.3 in \cite{huang2006sequential} and the discussions in \cite{scholtes2001convergence}.
    From Lemma \ref{lemma: feasibility of limit point}, we see that $\bar{z}$ is a feasible point of MPCC (\ref{equation: MPCC}).
    Since $z^k$ is a stationary point of $\mathcal{P}_{gap}(\mu^k)$, there exists a multiplier $u^k  \in \mathbb{R}^{n_h}$ such that
    \begin{equation}\label{equation: stationary condition of each iterate (original)}
        \nabla_z J(z^k) + (u^k)^T \nabla_z h(z^k) + \mu^k \nabla_z \varphi^{ab}(z^k) = 0.
    \end{equation}
    Since $\nabla_z \varphi^{ab}(z^k) = \nabla_{\lambda} \varphi^{ab}(z^k) \nabla_z \lambda^k + \nabla_{\eta} \varphi^{ab}(z^k) \nabla_z \eta^k$, by defining variables $v^k,w^k \in \mathbb{R}^{n_{\lambda}}$ with
    \begin{equation}\label{equation: definition of multiplier sequence}
         v^k = -\mu^k (\nabla_{\lambda} \varphi^{ab}(z^k))^T,  w^k = -\mu^k (\nabla_{\eta} \varphi^{ab}(z^k))^T,
    \end{equation}
    condition (\ref{equation: stationary condition of each iterate (original)}) becomes  
    \begin{equation}\label{equation: stationary condition of each iterate}
        \nabla_z J(z^k) + (u^k)^T \nabla_z h(z^k) - (v^k)^T \nabla_z \lambda^k - (w^k)^T \nabla_z \eta^k = 0.
    \end{equation}
    In the following, \textcolor{blue}{we prove the existence of multipliers satisfying (\ref{equation: Clarke stationarity analysis equation}) by showing that sequences $\{u^k\}_{k=0}^{\infty}$, $\{v^k\}_{k=0}^{\infty}$ and $\{w^k\}_{k=0}^{\infty}$ are bounded; moreover, their limit points satisfy the definition of Clarke stationarity (\ref{equation: Clarke stationarity analysis equation})}.
    
    The boundness of $\{u^k\}_{k=0}^{\infty}$, $\{v^k\}_{k=0}^{\infty}$ and $\{w^k\}_{k=0}^{\infty}$ is proved by contradiction.
    If these sequences are unbounded, we can find a subsequence $\mathcal{K}$ such that the normed sequence converges \cite{scholtes2001convergence}: 
    $\lim\limits_{k \rightarrow \infty} \frac{u^k}{\| \tau^k \|_1} = \bar{u}$,
    $\lim\limits_{k \rightarrow \infty} \frac{v^k}{\| \tau^k \|_1} = \bar{v}$, and
    $\lim\limits_{k \rightarrow \infty} \frac{w^k}{\| \tau^k \|_1} = \bar{w}$ for $k \in \mathcal{K}$,
    where $\tau = [u^T, v^T, w^T]^T$ and $\lim\limits_{k \rightarrow \infty} \| \tau^k \|_1 = + \infty$ for $ k \in \mathcal{K}$.
    Additionally, we have 
    \begin{equation}\label{equation: sum of multiplier limit}
        \| \bar{u} \|_1 + \| \bar{v} \|_1 + \| \bar{w} \|_1 = 1.
    \end{equation}
    Dividing (\ref{equation: stationary condition of each iterate}) by $\tau^k$ and taking the limit as $k \rightarrow \infty$, we have
    \begin{equation}\label{equation: stationary condition of each iterate (limit)}
        (\bar{u})^T \nabla_z h(\bar{z}) - (\bar{v})^T \nabla_z \bar{\lambda}  - (\bar{w})^T \nabla_z \bar{\eta}= 0.
    \end{equation}  
    Next, we show that \textcolor{blue}{the elements of limit points $\bar{v}$ and $\bar{w}$ satisfy} 
    \begin{equation}\label{equation: zero multiplier constraint inactived}
        \bar{v}_i = 0, \  i \in \mathcal{I}_{\eta}(\bar{z}) \text{ and } \bar{w}_i = 0, \  i \in \mathcal{I}_{\lambda}(\bar{z}).
    \end{equation}
    \textcolor{blue}{From (\ref{equation: reformulate D gap function as the sum of scalar functions}) and (\ref{equation: definition of multiplier sequence}), the elements of $v^k$ and $w^k$ are given by}
    \begin{equation}\label{equation: definition of multiplier sequence (elements)}
        v^k_i = -\mu^k \nabla_{\lambda_i} \delta^{ab}(\lambda^k_i, \eta^k_i), \ w^k_i = -\mu^k \nabla_{\eta_i} \delta^{ab}(\lambda^k_i, \eta^k_i).
    \end{equation} 
    \textcolor{blue}{Thus, the proof of (\ref{equation: zero multiplier constraint inactived}) requires analyzing $\nabla_{\lambda_i} \delta^{ab}(\lambda_i, \eta_i)$ and $\nabla_{\eta_i} \delta^{ab}(\lambda_i, \eta_i)$,
    which is expanded based on (\ref{equation: explicit expression of D gap gradient}) and (\ref{equation: scalar D gap function explicit expression}):}
    \begin{equation}\label{equation: expand scalar D gap grad wrt lambda}
        \nabla_{\lambda_i} \delta^{ab}(\lambda_i, \eta_i) =
            \begin{cases}
                (b-a) \lambda_i,      & \text{if } \eta_i \geq b \lambda_i, \  \eta_i \geq a \lambda_i, \\
                \eta_i - a \lambda_i, & \text{if } \eta_i < b \lambda_i, \  \eta_i > a \lambda_i, \\
                0,                    & \text{if } \eta_i \leq b \lambda_i, \  \eta_i \leq a \lambda_i, \\
                b \lambda_i - \eta_i, & \text{if } \eta_i > b \lambda_i, \  \eta_i < a \lambda_i, \\
            \end{cases}
    \end{equation}  
    \begin{equation}\label{equation: expand scalar D gap grad wrt eta}
         \nabla_{\eta_i} \delta^{ab}(\lambda_i, \eta_i) = 
            \begin{cases}
                0,                                 & \text{if } \eta_i \geq b \lambda_i, \  \eta_i \geq a \lambda_i, \\
                \frac{b\lambda_i - \eta_i}{b} ,    & \text{if } \eta_i < b \lambda_i, \  \eta_i > a \lambda_i, \\
                (\frac{1}{a} - \frac{1}{b})\eta_i, & \text{if } \eta_i \leq b \lambda_i, \  \eta_i \leq a \lambda_i, \\
                \frac{\eta_i - a \lambda_i}{a},    & \text{if } \eta_i > b \lambda_i, \  \eta_i < a \lambda_i. \\
            \end{cases}   
    \end{equation}    
    \textcolor{blue}{Then, we can determine the sign of $v^k_i$ and $w^k_i$ based on (\ref{equation: definition of multiplier sequence (elements)})--(\ref{equation: expand scalar D gap grad wrt eta}), $\mu^k > 0$, and $b > a > 0$, as illustrated in Fig. \ref{fig:sign of multiplier}}.
    According to the sign of $v^k_i$ and $w^k_i$, we can now state the proof of (\ref{equation: zero multiplier constraint inactived}) as follows.
    \begin{figure}
        \centering
        \includegraphics[width=1\linewidth]{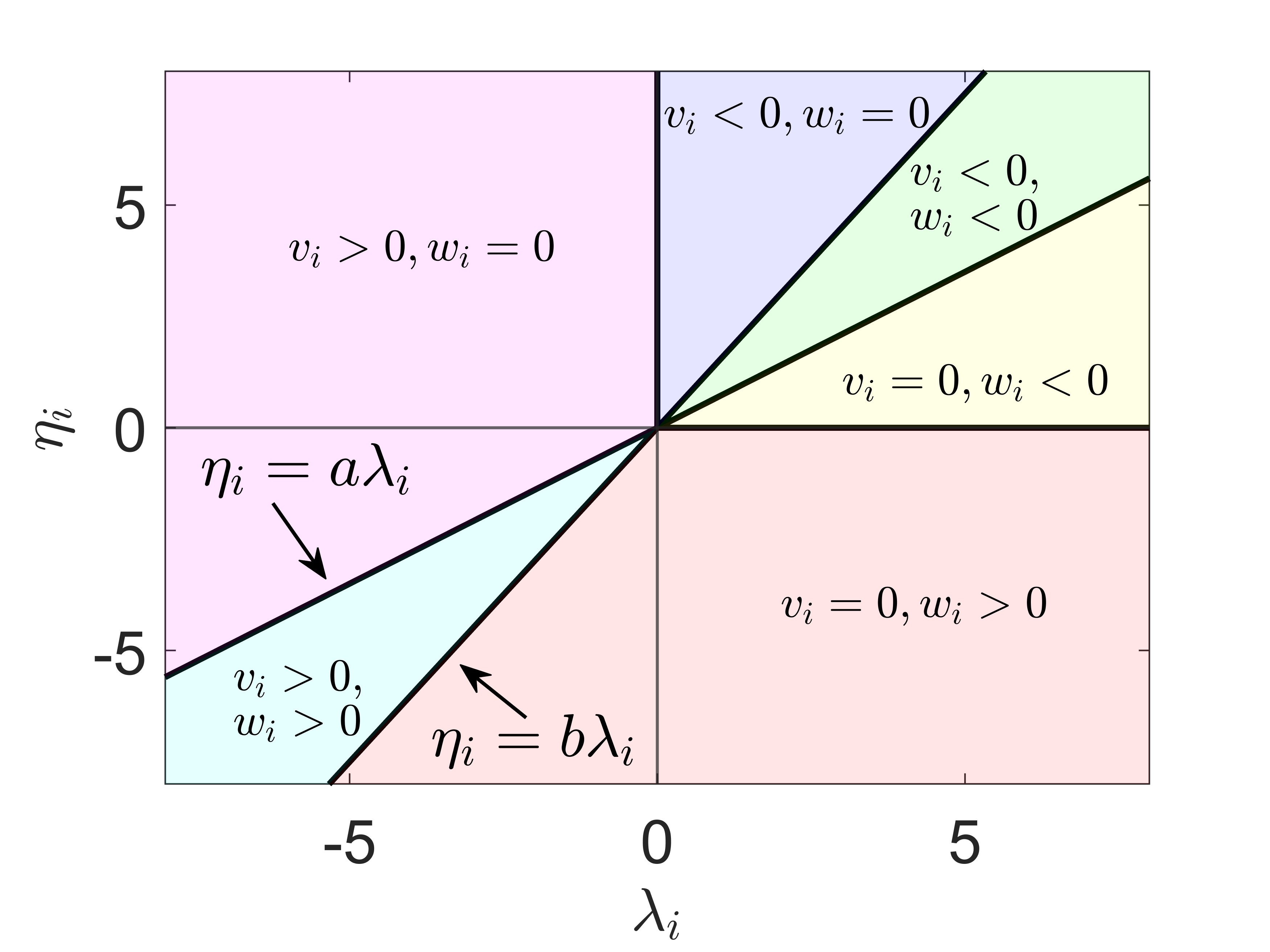}
        \caption{\textcolor{blue}{$\lambda_i$-$\eta_i$ plane is divided into six regions (distinguished by different colors) by axes $\lambda_i = 0$ (nonnegative part), $\eta_i = 0$ (nonnegative part), $\eta_i = a\lambda_i$, and $\eta_i = b\lambda_i$, where the signs of $v_i$ and $w_i$ vary in each region.}}
        \label{fig:sign of multiplier}
    \end{figure}

    \textcolor{blue}{Regarding $\bar{v}_i = 0,i \in \mathcal{I}_{\eta}(\bar{z})$, since $\bar{z}$ is a feasible point of MPCC (\ref{equation: MPCC}), the elements of $\bar{z}$ corresponding to $i \in \mathcal{I}_{\eta}(\bar{z})$ satisfy $\bar{\lambda}_i > 0, \bar{\eta}_i = 0$, and we have $\lambda^k_i \rightarrow \bar{\lambda}_i > 0$ and $\eta^k_i \rightarrow \bar{\eta}_i = 0$ for $i \in \mathcal{I}_{\eta}(\bar{z})$. 
    In other words, the pair $(\lambda^k_i, \eta^k_i)$ with $i \in \mathcal{I}_{\eta}(\bar{z})$ converges to the positive part of the axis $\eta_i = 0$. 
    Along with the sign of $v_i$ illustrated in Fig. \ref{fig:sign of multiplier}, we have $v^k_i = 0$ when $k$ is sufficiently large.}
    Consequently, we have
    \begin{equation}
        \lim\limits_{k \rightarrow \infty} \frac{v^k_i}{ \| \tau^k\|_1}  = \bar{v}_i = 0, \quad i \in \mathcal{I}_{\eta}(\bar{z}).
    \end{equation}
    
    \textcolor{blue}{Regarding $\bar{w}_i = 0, i \in \mathcal{I}_{\lambda}(\bar{z})$, Since $\bar{z}$ is a feasible point of MPCC (\ref{equation: MPCC}), the elements of $\bar{z}$ corresponding to $i \in \mathcal{I}_{\lambda}(\bar{z})$ satisfy  $\bar{\lambda}_i = 0, \bar{\eta}_i > 0$, and we have $\lambda^k_i \rightarrow \bar{\lambda}_i = 0$ and $\eta^k_i \rightarrow \bar{\eta}_i > 0$ for $i \in \mathcal{I}_{\lambda}(\bar{z})$.
    In other words, the pair $(\lambda^k_i, \eta^k_i)$ with $i \in \mathcal{I}_{\lambda}(\bar{z})$ converges to the positive part of the axis $\lambda_i = 0$.
    Along with the sign of $w_i$ illustrated in Fig. \ref{fig:sign of multiplier}, we have $w^k_i = 0$ when $k$ is sufficiently large.}
    Consequently, we have
    \begin{equation}
        \lim\limits_{k \rightarrow \infty} \frac{w^k_i}{ \| \tau^k\|_1} = \bar{w}_i = 0, \quad i \in \mathcal{I}_{\lambda}(\bar{z})
    \end{equation}

    Consequently, from (\ref{equation: sum of multiplier limit}) -- (\ref{equation: zero multiplier constraint inactived}) we have
    \begin{equation}
        \begin{split}
            &  \sum_{i = 1}^{n_h} \bar{u}_i \nabla_z h_i (\bar{z}) 
               + \sum_{i \in \mathcal{I}_{\lambda} (\bar{z}) \cup \mathcal{I}_{\lambda \eta}  (\bar{z})} \bar{v}_i \nabla_z \lambda_i \\ 
            & + \sum_{i \in \mathcal{I}_{\eta} (\bar{z}) \cup \mathcal{I}_{\lambda \eta} (\bar{z}) } \bar{w}_i \nabla_z \eta_i = 0,
        \end{split}       
    \end{equation}  
    with 
    \begin{equation}
        \sum_{i = 1}^{n_h} |\bar{u}_i| 
        + \sum_{i \in \mathcal{I}_{\lambda} (\bar{z}) \cup \mathcal{I}_{\lambda \eta}  (\bar{z})} | \bar{v}_i | 
        + \sum_{i \in \mathcal{I}_{\eta} (\bar{z}) \cup \mathcal{I}_{\lambda \eta} (\bar{z}) } | \bar{w}_i | = 1,
    \end{equation}   
    which contradicts the assumption that the MPCC-LICQ holds at $\bar{z}$.
    \textcolor{blue}{Thus, $\{u^k\}_{k=0}^{\infty}$, $\{v^k\}_{k=0}^{\infty}$ and $\{w^k\}_{k=0}^{\infty}$ are bounded and have limit points by the Bolzano–Weierstrass theorem. 
    Let the limit point of each sequence be
    $\bar{u} = \lim\limits_{k \rightarrow \infty} u^k$, 
    $\bar{v} = \lim\limits_{k \rightarrow \infty} v^k$, and 
    $\bar{w} = \lim\limits_{k \rightarrow \infty} w^k$ for $k \in \mathcal{K}$}.  
    Taking the limit of (\ref{equation: stationary condition of each iterate}) as $k \rightarrow \infty$ for $k \in \mathcal{K}$, we have
    \begin{equation}\label{equation: stationary condition of each iterate with multipliers (limit)}
        \nabla_z J(\bar{z}) + (\bar{u})^T \nabla_z h(\bar{z}) - (\bar{v})^T \nabla_z \bar{\lambda}  - (\bar{w})^T \nabla_z \bar{\eta}= 0,
    \end{equation}
    which is (\ref{equation: Clarke stationarity analysis equation stationary condition}). 
    Similar to the proof of (\ref{equation: zero multiplier constraint inactived}), we have
    \begin{equation}
        \lim_{k \rightarrow \infty} w^k_i = \bar{w}_i = 0, i \in \mathcal{I}_{\lambda}(\bar{z}), \ 
        \lim_{k \rightarrow \infty} v^k_i  = \bar{v}_i = 0, i \in \mathcal{I}_{\eta}(\bar{z}),
    \end{equation}
    which are (\ref{equation: Clarke stationarity analysis equation active lambda}) and (\ref{equation: Clarke stationarity analysis equation active eta}), respectively.    
    From (\ref{equation: expand scalar D gap grad wrt lambda}) and (\ref{equation: expand scalar D gap grad wrt eta}), we have
    \begin{equation*}
        \begin{split}
            0 \leq   & \nabla_{\lambda_i} \delta^{ab}(\lambda_i, \eta_i) \nabla_{\eta_i} \delta^{ab}(\lambda_i, \eta_i) \\
            = & 
            \begin{cases}
                0,                                                         & \text{if } \eta_i \geq b \lambda_i, \  \eta_i \geq a \lambda_i, \\
                (\eta_i - a \lambda_i)(\frac{b\lambda_i - \eta_i}{b}),    & \text{if } \eta_i < b \lambda_i, \  \eta_i > a \lambda_i, \\
                0,                                                         & \text{if } \eta_i \leq b \lambda_i, \  \eta_i \leq a \lambda_i, \\
                (b \lambda_i - \eta_i)(\frac{\eta_i - a \lambda_i}{a}),    & \text{if } \eta_i > b \lambda_i, \  \eta_i < a \lambda_i, \\
            \end{cases}
        \end{split}
    \end{equation*}
    This implies that for all $i \in \{1, \cdots, n_{\lambda} \}$, we have
    \begin{equation*}
        v^k_i w^k_i = (\mu^k)^2 \nabla_{\lambda_i} \delta^{ab}(\lambda^k_i, \eta^k_i) \nabla_{\eta_i} \delta^{ab}(\lambda^k_i, \eta^k_i)  \geq 0, 
    \end{equation*}
    and their limit are $\lim\limits_{k \rightarrow \infty} v^k_i w^k_i = \bar{v}_i \bar{w}_i \geq 0, k \in \mathcal{K}$.
    Thus (\ref{equation: Clarke stationarity analysis equation active lambda eta}) also holds at $\bar{z}$.
    This completes the proof.
\end{IEEEproof}

\subsection{Strong stationarity analysis}
Now, we present the main result of this study: $\bar{z}$ is a strongly stationary point of MPCC (\ref{equation: MPCC}) if the additional assumptions hold, that is, $z^k$ satisfies the second-order necessary condition of $\mathcal{P}_{gap}(\mu^k)$ and the ULSC condition holds at $\bar{z}$.

\begin{theorem}\label{theorem: strong stationarity of limit point}
    For a given $\mu^k > 0$, let $z^k$ be a stationary point of $\mathcal{P}_{gap}(\mu^k)$, and assume that $z^k$ satisfies the second-order necessary condition of $\mathcal{P}_{gap}(\mu^k)$.
    Let $\mu^k \rightarrow + \infty$ and $\bar{z}$ be a limit point of $\{z^k\}_{k = 0}^{\infty}$.    
    Let the assumption of Theorem \ref{theorem: Clarke stationarity of limit point} hold, and assume that the ULSC condition holds at $\bar{z}$.    
    Then, $\bar{z}$ is a strongly stationary point of MPCC (\ref{equation: MPCC}), that is, there exist Lagrangian multipliers $\bar{u} \in\mathbb{R}^{n_h}$, $\bar{v},\bar{w} \in \mathbb{R}^{n_{\lambda}}$ such that:
        \begin{subequations}\label{equation: strong stationarity analysis equation}
            \begin{align}
                & \nabla_z J(\bar{z}) + \bar{u}^T \nabla_z h(\bar{z}) - \bar{v}^T \nabla_z \bar{\lambda} - \bar{w}^T \nabla_z \bar{\eta} = 0, \label{equation: strong stationarity analysis equation stationary condition}\\
                & \bar{v}_i \in \mathbb{R}, \bar{w}_i = 0, \ i \in \mathcal{I}_{\lambda}(\bar{z}), \label{equation: strong stationarity analysis equation active lambda}\\
                & \bar{v}_i = 0,            \bar{w}_i \in \mathbb{R},  \ i \in \mathcal{I}_{\eta}(\bar{z}), \label{equation: strong stationarity analysis equation active eta}\\
                & \bar{v}_i \geq 0,         \bar{w}_i \geq 0,  \ i \in \mathcal{I}_{\lambda \eta}(\bar{z}) \label{equation: strong stationarity analysis equation active lambda eta}.
            \end{align}
        \end{subequations}    
\end{theorem}
\begin{IEEEproof}
    We only need to prove that (\ref{equation: strong stationarity analysis equation active lambda eta}) holds at $\bar{z}$.

    From Theorem \ref{theorem: Clarke stationarity of limit point}, we have $\bar{v}_i \bar{w}_i \geq 0, \ i \in \mathcal{I}_{\lambda \eta}(\bar{z})$.
    Since the ULSC condition holds at $\bar{z}$, that is, $\bar{v}_i \bar{w}_i \neq 0, \forall i \in \mathcal{I}_{\lambda \eta}(\bar{z})$, we have that $\bar{v}_i \bar{w}_i \geq 0, i \in \mathcal{I}_{\lambda \eta}(\bar{z})$ holds if either $\bar{v}_i >0, \bar{w}_i > 0$ or $\bar{v}_i < 0, \bar{w}_i < 0, i \in \mathcal{I}_{\lambda \eta}(\bar{z})$.
    In the following, we prove by contradiction that it can only be $\bar{v}_i >0, \bar{w}_i > 0$.
    
    First, if $\bar{v}_i >0, \bar{w}_i > 0$,  then when $k$ is sufficiently large, we have $v^k_i >0, w^k_i > 0$, 
    moreover, following from the sign of $v_i$ and $w_i$ in Fig. \ref{fig:sign of multiplier}, we have $\eta^k_i > b \lambda^k_i$ and $\eta^k_i < a \lambda^k_i$.
    Similarly, if $\bar{v}_i < 0, \bar{w}_i < 0$, we have $\eta^k_i < b \lambda^k_i$ and $\eta^k_i > a \lambda^k_i$ when $k$ is sufficiently large.
    Thus, we can transform the analysis of the multiplier sequence $\{v^k_i\}, \{w^k_i\}$ into the analysis of the primal variable sequences $\{\lambda^k_i\}, \{\eta^k_i\}$.
 
    Next, we investigate the primal variable sequences.
    Since $z^k$ satisfies the second-order necessary condition of $\mathcal{P}_{gap}(\mu^k)$, for any $d \in \mathcal{C}(z^k) := \{ d \in \mathbb{R}^{n_x + 2 n_{\lambda}} | \nabla_z h(z^k) d = 0\}$, we have:
    \begin{equation}\label{equation: gap penalty reformulation second order necessary optimality condition}
        \begin{split}
             0 \leq & D^2(J(z^k); d, d) + D^2(\mu^k\varphi^{ab}(z^k); d, d) \\
             = & D^2(\mu^k\varphi^{ab}(z^k); d, d) - D(\mu^k\varphi^{ab}(z^k); d) \\ 
             & + \frac{1}{2}d^T \nabla_{zz}J(z^k)d\\
             = & \mu^k\sum^{n_{\lambda}}_{i = 1} \{ D^2(\delta^{ab}(\lambda_{i}^k, \eta_{i}^k); d, d) - D(\delta^{ab}(\lambda_{i}^k, \eta_{i}^k); d) \} \\
             & + \frac{1}{2}d^T \nabla_{zz}J(z^k)d 
        \end{split}    
    \end{equation}
    Let $d_{\lambda_i}$ and $d_{\eta_i}$ be the elements of $d$ associated with $\lambda_i$ and $\eta_i$, respectively.
    Following from Proposition \ref{proposition: properties of the directional derivatives}, we have
    \begin{equation}\label{equation: the difference between the first and second order directional derivative}
        \begin{split}
            & D^2(\delta^{ab}(\lambda_{i}, \eta_{i}); d, d) - D(\delta^{ab}(\lambda_{i}, \eta_{i}); d) \\
            = & 
            \begin{cases}
               \frac{b-a}{2}d_{\lambda_i}^2 & \text{if } \eta_i > b \lambda_i,  \eta_i > a \lambda_i \\
               \frac{b-a}{2}d^2_{\lambda_i} - \frac{1}{2b}(d_{\eta_i} - b d_{\lambda_i})^2 + m_b     & \text{if } \eta_i = b \lambda_i,  \eta_i > a \lambda_i \\
               \frac{b-a}{2}d^2_{\lambda_i} + \frac{1}{2a}(d_{\eta_i} - a d_{\lambda_i})^2 - m_a     & \text{if } \eta_i > b \lambda_i,  \eta_i = a \lambda_i \\
               - \frac{a}{2}d^2_{\lambda_i} + d_{\lambda_i}d_{\eta_i} - \frac{1}{2b}d^2_{\eta_i}    & \text{if } \eta_i < b \lambda_i,  \eta_i > a \lambda_i \\
               \frac{b-a}{2ab}d_{\eta_i}^2    & \text{if } \eta_i < b \lambda_i,  \eta_i < a \lambda_i \\
               \frac{b-a}{2ab}d_{\eta_i}^2 + m_b    & \text{if } \eta_i = b \lambda_i,  \eta_i < a \lambda_i \\
               \frac{b-a}{2ab}d_{\eta_i}^2 - m_a     & \text{if } \eta_i < b \lambda_i,  \eta_i = a \lambda_i \\
               \frac{b}{2}d^2_{\lambda_i} - d_{\lambda_i}d_{\eta_i} + \frac{1}{2a}d^2_{\eta_i}   & \text{if } \eta_i > b \lambda_i,  \eta_i < a \lambda_i                
            \end{cases}
        \end{split}
    \end{equation}
    where $m_a, m_b: \mathbb{R} \times \mathbb{R} \rightarrow \mathbb{R}$ are functions defined by 
    $m_a(d_{\lambda_i}, d_{\eta_i}) = \frac{1}{2a}(\max(0,d_{\eta_i} - a d_{\lambda_i}))^2$ and 
    $m_b(d_{\lambda_i}, d_{\eta_i}) = \frac{1}{2b}(\max(0,d_{\eta_i} - b d_{\lambda_i}))^2$, respectively.

    Here, we claim that the elements of the sequence $\{\lambda^k_i\}, \{\eta^k_i\}$ can neither satisfy $\eta^k_i = a \lambda^k_i$ nor $\eta^k_i = b \lambda^k_i$ when $k$ is sufficiently large.
    Otherwise, if $\lambda^k_i$ and $\eta^k_i$ satisfy either $\eta^k_i = a \lambda^k_i$ or $\eta^k_i = b \lambda^k_i$ when $k$ is sufficiently large, we have $\lambda^k_i, \eta^k_i \rightarrow 0$, which indicates that $i \in \mathcal{I}_{\lambda \eta}(\bar{z})$.
    However, following from the definition and sign of $v_i$ and $w_i$ in (\ref{equation: definition of multiplier sequence (elements)}) and Fig. \ref{fig:sign of multiplier}, the associated multipliers $v^k_i$ or $w^k_i$ will equal zero when $k$ is sufficiently large, that is
    \begin{equation*}
        \begin{cases}
            w^k_i = 0, & \text{if } \eta^k_i = b \lambda^k_i,  \eta^k_i > a \lambda^k_i \text{ or } \eta^k_i > b \lambda^k_i,  \eta^k_i = a \lambda^k_i\\
            v^k_i = 0, & \text{if } \eta^k_i = b \lambda^k_i,  \eta^k_i < a \lambda^k_i \text{ or } \eta^k_i < b \lambda^k_i,  \eta^k_i = a \lambda^k_i 
        \end{cases}
    \end{equation*}
    and we have $v^k_i \rightarrow \bar{v}_i = 0$ or $w^k_i \rightarrow \bar{w}_i = 0$. 
    This contradicts the assumption that the ULSC condition holds at $\bar{z}$.
    For brevity, we define the following index set
    \begin{subequations}
        \begin{align}
            & \mathcal{J}_1(z) = \{ {i \in \{1, \cdots, n_{\lambda} \} | \eta_i > b \lambda_i,  \eta_i > a \lambda_i\}}, \\
            & \mathcal{J}_2(z) = \{ {i \in \{1, \cdots, n_{\lambda} \} | \eta_i < b \lambda_i,  \eta_i > a \lambda_i\}}, \\
            & \mathcal{J}_3(z) = \{ {i \in \{1, \cdots, n_{\lambda} \} | \eta_i < b \lambda_i,  \eta_i < a \lambda_i\}}, \\
            & \mathcal{J}_4(z) = \{ {i \in \{1, \cdots, n_{\lambda} \} | \eta_i > b \lambda_i,  \eta_i < a \lambda_i\}},
        \end{align}
    \end{subequations}    
    and have that $ \mathcal{I}_{\lambda}(\bar{z}) = \mathcal{J}_1(z^k)$ and $\mathcal{I}_{\eta}(\bar{z}) = \mathcal{J}_3(z^k)$ when $k$ is sufficiently large.
    Assuming that $\bar{v}_i < 0, \bar{w}_i < 0, i \in \mathcal{I}_{\lambda \eta}(\bar{z})$, we have $\mathcal{I}_{\lambda \eta}(\bar{z}) = \mathcal{J}_2(z^k)$ when $k$ is sufficiently large.
    Substituting (\ref{equation: the difference between the first and second order directional derivative}) into (\ref{equation: gap penalty reformulation second order necessary optimality condition}), we have
    \begin{equation}\label{equation: gap penalty reformulation second order necessary optimality condition impossible case}
        \begin{split}
            0 \leq \ & \mu^k  \sum_{i \in \mathcal{J}_1(z^k)} (\frac{b-a}{2}d_{\lambda_i}^2) + \mu^k \sum_{i \in \mathcal{J}_3(z^k)}( \frac{b-a}{2ab}d_{\eta_i}^2)  \\ 
                    & + \mu^k \sum_{i \in \mathcal{J}_2(z^k)} (- \frac{a}{2}d^2_{\lambda_i} + d_{\lambda_i}d_{\eta_i} - \frac{1}{2b}d^2_{\eta_i})  \\
                    & + \frac{1}{2}d^T \nabla_{zz}J(z^k)d 
        \end{split}
    \end{equation}
    Since the MPCC-LICQ holds at $\bar{z}$, there exists $d$ such that
    \begin{subequations}\label{equation: choice of d}
        \begin{align}
            & \nabla_z h(z^k) d = 0, \\
            & \nabla_z \lambda^k_i d = d_{\lambda_i} = 0, \ i \in \mathcal{I}_{\lambda} (\bar{z}), \\
            & \nabla_z \eta^k_i d = d_{\eta_i} = 0, \ i \in \mathcal{I}_{\eta} (\bar{z}),  \\
            & \nabla_z \lambda^k_i d = d_{\lambda_i} > 0, \ i \in \mathcal{I}_{\lambda \eta}(\bar{z}), \\
            & \nabla_z \eta^k_i d = d_{\eta_i} < 0, \ i \in \mathcal{I}_{\lambda\eta} (\bar{z}), 
        \end{align}
    \end{subequations}
    when $k$ is sufficiently large.
    However, since we have
    \begin{equation*}
        - \frac{a}{2}d^2_{\lambda_i} + d_{\lambda_i}d_{\eta_i} - \frac{1}{2b}d^2_{\eta_i} < 0,  \ \forall d_{\lambda_i}d_{\eta_i}<0,
    \end{equation*} 
    the second-order necessary condition (\ref{equation: gap penalty reformulation second order necessary optimality condition impossible case}) fails to hold when $k$ is sufficiently large.
    This contradicts the assumption that $z^k$ satisfies the second-order necessary condition of $\mathcal{P}_{gap}(\mu^k)$.
    
    In fact, if assume $\bar{v}_i > 0, \bar{w}_i > 0, i \in \mathcal{I}_{\lambda \eta}(\bar{z})$, then $\mathcal{I}_{\lambda \eta}(\bar{z}) = \mathcal{J}_4(z^k)$ when $k$ is sufficiently large.
    Similarly, (\ref{equation: gap penalty reformulation second order necessary optimality condition}) becomes
    \begin{equation*}
        \begin{split}
            0 \leq \ & \mu^k \sum_{i \in \mathcal{J}_4(z^k)} (\frac{b}{2}d^2_{\lambda_i} - d_{\lambda_i}d_{\eta_i} + \frac{1}{2a}d^2_{\eta_i}) + \frac{1}{2}d^T \nabla_{zz}J(z^k)d                
        \end{split}    
    \end{equation*}
    which holds even when $k$ is sufficiently large because
    \begin{equation*}
        \frac{b}{2}d^2_{\lambda_i} - d_{\lambda_i}d_{\eta_i} + \frac{1}{2a}d^2_{\eta_i} \geq 0, \ \forall d_{\lambda_i}, d_{\eta_i} \in \mathbb{R}.
    \end{equation*}
    This completes the proof.  
\end{IEEEproof}

\section{Conclusion}\label{section: conclusion}
This study provides a theoretical guarantee for the convergence properties of our previous work \cite{lin2024gap}, where a new penalty method is proposed to solve MPCC. 
We proved that the new penalty method converges to a Clarke stationary point of MPCC, provided that the MPCC-LICQ holds.
It converges to a strongly stationary point of MPCC, provided that, additionally, the ULSC condition holds and the gap penalty reformulation satisfies the second-order necessary condition in terms of the second-order directional derivatives.
Convergence to a strongly stationary point of MPCC is favorable as it indicates that our new penalty method can find an MPCC solution.
\textcolor{blue}{Future work will consider an MPCC with inequality constraints and focus on the convergence analysis to the MPCC second-order optimality conditions}.

\appendix

\subsection{Proof of the second statement in Theorem \ref{theorem: necessary optimality conditions of the NLP problem with LC1 function}}\label{subsection: proof of necessary optimality conditions of the NLP problem with LC1 function}

\begin{IEEEproof}
    Let $x_{t,d} = x^* + td + t^2d$ be a point near $x^*$ with $t \geq 0$, and $d_{x, \gamma} = [d^T, 0_{1 \times n_h} ]^T$ be the direction with $d \in \mathcal{C}(x^*)$.
    First, we use (\ref{equation: function expansion in terms of directional derivatives}) to expand the Lagrangian $\mathcal{L}(x_{t,d}, \gamma^*)$:
    \begin{equation}\label{equation: Lagrangian expansion in terms of directional derivatives}
        \begin{split}
            \mathcal{L}(x_{t,d}, \gamma^*) = & \mathcal{L}(x^*, \gamma^*) + t D(\mathcal{L}(x^*, \gamma^*); d_{x, \gamma}) \\
                                             &  + t^2 D^2(\mathcal{L}(x^*, \gamma^*); d_{x, \gamma}, d_{x, \gamma}) + o(t^2).
        \end{split}
    \end{equation}
    Regarding the term on the left-hand-side of (\ref{equation: Lagrangian expansion in terms of directional derivatives}), 
    \begin{equation*}\label{equation: Lagrangian expansion lhs term}
        \begin{split}
            & \mathcal{L}(x_{t,d}, \gamma^*) \\ 
            = & J(x_{t,d}) + g(x_{t,d}) + \sum_{i = 1}^{m} \gamma^*_i h_i(x_{t,d}) \\
            = & J(x_{t,d}) + g(x_{t,d}) + \sum_{i = 1}^{m} \gamma^*_i \{ \underbrace{h_i(x^*)}_{= 0} + t \underbrace{D(h_i(x^*); d)}_{= \nabla_xh_i(x^*) d = 0} \\ 
              &  + t^2 D^2(h_i(x^*); d, d) + o(t^2)\} \\
            = & J(x_{t,d}) + g(x_{t,d}) + t^2 \sum_{i = 1}^{m} \gamma^*_i  D^2(h_i(x^*); d, d) + o(t^2).
        \end{split}
    \end{equation*}
    Regarding the first term on the right-hand-side of (\ref{equation: Lagrangian expansion in terms of directional derivatives}),
    \begin{equation*}\label{equation: Lagrangian expansion rhs term 1}
        \mathcal{L}(x^*, \gamma^*) = J(x^*) + g(x^*) + \sum_{i = 1}^{m} \gamma^*_i h_i(x^*) = J(x^*) + g(x^*).
    \end{equation*}
    Regarding the second term on the right-hand-side of (\ref{equation: Lagrangian expansion in terms of directional derivatives}), 
    \begin{equation*}\label{equation: Lagrangian expansion rhs term 2}
        t D(\mathcal{L}(x^*, \gamma^*); d_{x, \gamma}) = t \underbrace{\nabla_x\mathcal{L}(x^*, \gamma^*)}_{= 0} d = 0.
    \end{equation*}
    Regarding the third term on the right-hand-side of (\ref{equation: Lagrangian expansion in terms of directional derivatives}), 
    \begin{equation*}\label{equation: Lagrangian expansion rhs term 3}
        \begin{split}
              & t^2 D^2(\mathcal{L}(x^*, \gamma^*); d_{x, \gamma}, d_{x, \gamma}) \\
            = & t^2 D^2(J(x^*) + g(x^*); d, d) + t^2 \sum_{i = 1}^{m} \gamma^*_i D^2(h_i(x^*); d, d).
        \end{split}
    \end{equation*}   
    Therefore, (\ref{equation: Lagrangian expansion in terms of directional derivatives}) can be rewritten as
    \begin{equation}\label{equation: Lagrangian expansion in terms of directional derivatives (simplified)}
        \begin{split}
             & J(x_{t,d}) + g(x_{t,d}) \\ 
             = & J(x^*) + g(x^*) + t^2 D^2(J(x^*) + g(x^*); d, d) + o(t^2).
        \end{split}        
    \end{equation}
    If $D^2(J(x^*) + g(x^*); d, d) < 0$, then (\ref{equation: Lagrangian expansion in terms of directional derivatives (simplified)}) would imply that $J(x_{t,d}) + g(x_{t,d}) < J(x^*) + g(x^*)$ for a sufficiently small $t$, which contradicts the fact that $x^*$ is a local minimum of (\ref{equation: NLP problem with LC1 function}). 
    Thus (\ref{equation: NLP problem with LC1 function second order necessary optimality condition}) is proven.    
\end{IEEEproof}


\bibliographystyle{IEEEtran}
\bibliography{IEEEabrv, reference}

\end{document}